\documentclass[english,12pt]{smfart}
  \usepackage{macros_anglais}
  \usepackage{dsfont}

\title{Approximation by finitely supported measures}
\author{Beno\^{\i}t Kloeckner}

\newcommand{\equivalent}{\mathop {\sim}\limits}
\newcommand{\supp}{\mathop {\mathrm{supp}}\nolimits}
\newcommand{\diam}{\mathop {\mathrm{diam}}\nolimits}
\newcommand{\vol}{\mathop {\mathrm{vol}}\nolimits}
\newcommand{\intd}{\mathrm{d}}

\begin{document}

\maketitle

\section{Introduction}

In this paper we are interested in the following question: given a finite
measure $\mu$, at what speed can it be approximated
by finitely supported measures? To give a sense to the question, one needs
a distance on the space of measures; we shall use the Wasserstein distances $W_p$,
with arbitrary exponent $p\in[1,+\infty)$
(definitions are recalled in Section \ref{sec:recalls}).

This problem has been called {\it Quantization for probability distribution},
the case of exponent $p=1$ has also been studied under the name of 
{\it location problem}, and
the case $p=2$ is linked with {\it optimal centroidal Voronoi tessellations}.
After submission of the present article, we became aware that the previous works
cover much more of the material presented than we first thought; see
Subsection \ref{sec:biblio} for detailled references.

This problem could be of interest for a numerical study of transportation problems,
where measures can be represented by discrete ones. One would need to know the number of
points needed to achieve some precision in the approximation. 

We shall restrict our attention to compactly supported Borelian measures on 
Riemannian manifolds.

\subsection{Statement of the results}

First we show that the order of convergence is determined by the dimension
of the measure (see definitions in Section \ref{sec:recalls}); 
$\Delta_N$ denotes the set of all measures supported in at most $N$ points.
\begin{theo}\label{theo:dimension}
If $\mu$ is compactly supported and Alhfors regular of dimension $s>0$, then
\[W_p(\mu,\Delta_N) \approx \frac1{N^{1/s}}.\]
\end{theo}
Here we write $\approx$ to say that one quantity is bounded above and below by 
positive multiples of the other. Examples of Ahlfors regular measures are given
by the volume measures on submanifolds, and self-similar measures (see for example \cite{Falconer}).
Theorem \ref{theo:dimension}, 
to be proved in a slightly more general and precise form in section \ref{sec:dimension},
is simple and unsurprising; it reminds of ball packing and covering, and indeed
relies on a standard covering argument.

In the particular case of absolutely continuous measures, one can give much finer estimates.
First, it is easily seen that if $\square^d$ denotes the uniform measure on a 
Euclidean unit cube of dimension $d$, then there is a constant $\theta(d,p)$ such that
\[W_p(\square^d,\Delta_N)\sim \frac{\theta(d,p)}{N^{1/d}}\]
(Proposition \ref{prop:cube}).
Note that determining the precise value of $\theta(d,p)$ seems difficult; known cases are discussed in Section 
\ref{sec:biblio}.

The main result of this paper is the following, where ``$\vol$'' denotes
the volume measure on the considered Riemannian manifold and is the default
measure for all integrals.
\begin{theo}\label{theo:cont}
If $\mu=\rho\vol$ where $\rho$ is a compactly supported function on a Riemannian manifold
$(M,g)$, then for all $1\leqslant p<\infty$ we have
\begin{equation}
W_p(\mu,\Delta_N) \equivalent \frac{\theta(d,p)\, |\rho|^{1/p}_{\frac d{d+p}}}{N^{1/d}}
\label{eq:main}
\end{equation}
where $|\rho|_\beta=(\int_M \rho^\beta)^{1/\beta}$ is the $L^\beta$ ``norm'',
here with $\beta<1$ though.

Moreover, if $(\mu_N)$ is a sequence of finitely supported measures such that $\mu_N\in\Delta_N$
minimizes $W_p(\mu,\mu_N)$, then the sequence of probability measures $(\bar\mu_N)$ that are 
uniform on the support of $\mu_N$ converges weakly to the multiple of 
$\rho^{\frac{d}{p+d}}$ that has mass $1$.
\end{theo}

Theorem \ref{theo:cont} is proved in Section \ref{sec:cont}. Note that the hypothesis 
that $\mu$ has compact support
is obviously needed: otherwise, $|\rho|_{d/(d+p)}$ could be infinite. 
Even when $\mu$ is in $L^{d/(d+p)}$,
there is the case where it is supported
on a sequence of small balls going to infinity. Then the location of the balls is 
important in the quality
of approximation and not only the profile of the density function. However,
this hypothesis could probably be relaxed to a moment condition.

Theorem \ref{theo:cont} has no real analog for measures
of fractional dimension.
\begin{theo}\label{theo:example}
There is a $s$-dimensional Ahlfors regular measure $\kappa$ on $\mathbb{R}$ (namely,
$\kappa$ is the Cantor dyadic measure) such that
$W_p(\kappa,\Delta_N) N^{1/s}$ has no limit.
\end{theo}
Section \ref{sec:examples} is devoted to this example.

Part of the interest of Theorem \ref{theo:cont} comes from the following observation, to be
discussed in Section \ref{sec:CVT}:
when $p=2$, the support of a distance minimizing $\mu_N\in\Delta_N$ generates a centroidal 
Voronoi tessellation,
that is, each point is the center of mass (with respect to $\mu$)
of its Voronoi cell. We thus get the asymptotic repartition
of an important family of centroidal Voronoi tessellations, which enables us to 
prove some sort of energy equidistribution principle.

\subsection{Discussion of previously known results}\label{sec:biblio}
There are several previous works closely related to the content of this paper.

\subsubsection{Foundations of Quantization for Probability Distributions}

The book \cite{Graf-Luschgy} by Graf and Luschgy (see also the references
therein), that we only discovered recently,
contains many results on the present problem.
Theorem \ref{theo:dimension} is proved there in section 12, but our proof
seems more direct. Theorem \ref{theo:cont} is proved in the Euclidean
case in Sections 6 and 7 (with a weakening of the compact support assumption).
A generalization of Theorem \ref{theo:example} is proved in Section 14, yet
we present a proof for the sake of self-completeness.

The case $p=1$, $M=\mathbb{R}^n$ is usually called the {\it location problem}.
In this setting, Theorem \ref{theo:cont}
has also been proved by Bouchitt\'e, Jimenez and Rajesh
\cite{Bouchitte} under the additionnal assumption that $\rho$ is lower semi-continuous.

Our main motivation to publish this work despite these overlaps
is that the case of measures on manifold
should find applications; for example, good approximations of the curvature
measure of a convex body by discrete measures should give good approximations
of the body by polyhedra.

It seems that the quantization, the
location problem and the study of optimal CVTs, although the last two
are particular cases of the first one, have been only studied independently. 
We hope that
notincing this proximity will encourage progress on each question to be 
translated in the others.

\subsubsection{Around the main theorem}

 Mosconi and
Tilli in \cite{Mosconi-Tilli} have studied (for any exponent $p$, in
$\mathbb{R}^n$) the {\it irrigation problem}, where the approximating measures
are supported on connected sets of length $<\ell$ (the length being the $1$-dimensional 
Hausdorff measure)
instead of being supported on $N$ points; the order of approximation is then $\ell^{1/(d-1)}$.

Brancolini, Buttazzo, Santambrogio and Stepanov compare in \cite{Brancolini} the location problem
with its ``short-term planning'' version, where the support of $\supp\mu_N$ is constructed
by adding one point to that of $\mu_{N-1}$, minimizing the cost only locally in $N$.

\subsubsection{Approximation constants for cubes}

Some values of $\theta(d,p)$ have been determined.
First, it is easy to compute them in dimension $1$:
\[\theta(1,p)=\frac{(p+1)^{-1/p}}{2}.\]
The case $d=2$ has been solved by Fejes T\'oth \cite{FejesToth,FejesToth2},
(and by Newmann \cite{Newman} for $p=2$ and Morgan and Bolton \cite{Morgan-Bolton}
for $p=1$), see also \cite{Graf-Luschgy} Section 8.
In particular
\[\theta(2,2)=5\sqrt{3}/54 \qquad \theta(2,1)=2^{-2/3}3^{-7/4}(4+\ln 27).\]
When $d=2$ and for all $p$, the hexagonal lattice is optimal
(that is, the given $\theta$ is the distance
between the uniform measure on a regular hexagon and a Dirac mass at its center).
All other cases are open to our knowledge.
For numerical evidence in the
case $p=2,d=3$ see Du and Wang \cite{Du-Wang}. Note that in the limit case
$p=\infty$, determining $\theta$ amounts to determining the minimal density of a ball covering of $\mathbb{R}^d$,
which is arguably as challenging as determining the maximal density of a ball packing, a well-known
open problem if $d>3$.

\subsubsection{Random variations}

Concerning the order of convergence, it is worth
comparing with the problem of estimating the distance from a measure $\mu$
to empirical measures $\bar\mu_N=N^{-1}\sum_k \delta_{X_k}$ where $X_1,\ldots, X_N$
are independent random variables of law $\mu$. It seems that $\bar\mu_N$ is almost optimal
in the sense that $W_2(\mu,\bar\mu_N)\sim C\, N^{-1/d}$ almost surely (under moment conditions, but here
we take $\mu$ compactly supported so this is not an issue); Horowitz and Karandikar have
shown in \cite{Horowitz-Karandikar} that $W_2(\mu,\bar\mu_N)$
has the order at most $N^{-1/(d+4)}$ and the better exponent above is suggested in the 
Mathematical Review of that paper. 

Let us also briefly note that the optimal matching problem for random data is
related to our problem. Simply put, one can say that if $\bar\mu'_N$ is another empirical measure
of $\mu$, then $W_2(\bar\mu_N,\bar\mu'_N)$ also has the order $N^{-1/d}$ if $d\geqslant 3$
(see for example Dobri\'c and Yukich \cite{Dobric-Yukich}). In the same
flavour, other optimisation problems for random data have been studied
(minimal length covering tree, traveling salesperson problem, bipartite version of those, etc.)

\subsubsection{Centroidal Voronoi Tesselations}

In the case $p=2$, the problem is linked to (optimal) {\it centroidal Voronoi Tesselation},
see Section \ref{sec:CVT} and \cite{Du-Faber-Gunzburger}.
In that paper (Section 6.4.1), the principle of energy equidistibution is given in 
the $1$-dimensional case for smooth density $\rho$. Our corollary \ref{coro:equid} 
in the last section generalize this to non regular densities, all exponents, and 
all dimensions; it is however quite a direct consequence
of Theorem \ref{theo:cont}.

\subsection{Related open questions}
The number $N$ of points of the support may be the first measure of complexity of a finitely
supported measure that one comes up with, but it is not necessarily the most relevant.
Concerning the problem of numerical analysis of transportation problems, numbers are usually
encoded in a computer by floating numbers. One could therefore define the complexity of a measure
supported on points of decimal coordinates, with decimal quantity of mass at each point as the
memory size needed to describe it, and search to minimize the distance to a given $\mu$ among
measures of given complexity.

Another possible notion of complexity is entropy : one defines
\[h\left(\sum_i m_i \delta_{x_i}\right)=-\sum_i m_i\ln(m_i).\]
A natural question is to search a $\mu_h$ that minimizes the distance to $\mu$ among the finitely 
supported measures of entropy at most $h$, and to study the behavior of $\mu_h$ when we let
$h\to\infty$.

\subsection*{Acknowledgements} I am grateful to Romain Joly, Vincent Mun\-nier,
Herv\'e Pajot, R\'emy Peyre and C\'edric Villani for interesting discussions or comments.

\section{Recalls and definitions}\label{sec:recalls}

\subsection{Notations}

Given two sequences $(u_n)$, $(v_n)$ of non-negative real numbers, we shall write:
\begin{itemize}
\item $u_n\lesssim v_n$ to mean that there exists a positive real $a$ and an integer $N_0$ such that
      $u_n \leqslant a v_n$ for all $n\geqslant N_0$,
\item $u_n\approx v_n$ to mean $u_n\lesssim v_n$ and $u_n\gtrsim v_n$.
\end{itemize}

From now on, $M$ is a given Riemannian manifold of dimension $d$. By a \emph{domain}
of $M$ we mean a compact domain with piecewise smooth boundary (and possibly corners)
and finitely many connected components.

\subsection{Ahlfors regularity and a covering result}

We denote by $B(x,r)$ the closed ball of radius $r$ and center $x$; sometimes,
when $B=B(x,r)$ and $k\in\mathbb{R}$, we denote by $kB$ the ball $B(x,kr)$.

Let $\mu$ be a finite, compactly supported measure on a manifold $M$ of
dimension $d$, and let $s\in(0,+\infty)$.
One says that $\mu$ is \emph{Ahlfors regular} of dimension $s$ if there is a constant $C$
such that for all $x\in\supp\mu$ and for all $r\leqslant\diam(\supp\mu)$,
one has 
\[C^{-1} r^s \leqslant \mu(B(x,r)) \leqslant C r^s.\]
This is a strong condition, but 
is satisfied for example by auto-similar measures, see \cite{Hutchinson,Falconer}
for definitions and Section \ref{sec:examples} for the most famous example of the Cantor measure.

Note that if $\mu$ is Ahlfors regular of dimension $s$,
then $s$ is the Hausdorff dimension of $\supp\mu$ (and therefore $s\leqslant d$),
see \cite[Sec. 8.7]{Heinonen}.

We shall need the following classical
covering result.
\begin{prop}[$5\delta$ covering]
If $X$ is a closed set and $\mathscr{F}$ is a family of balls of uniformly bounded diameter such that
$X\subset \bigcup_\mathscr{F} B$, then there is a subfamilly $\mathscr{G}$ of $\mathscr{F}$
such that:
\begin{itemize}
\item $X\subset \bigcup_\mathscr{G} 5B$,
\item $B\cap B'=\varnothing$ whenever $B\neq B' \in\mathscr{G}$.
\end{itemize}
\end{prop}

\subsection{Wasserstein distances}

Here we recall some basic facts on optimal transportation and Wasserstein distances.
For more information, the reader is suggested to look for example at Villani's book \cite{Villani}
which provides a very good introduction to this topic.

First consider the case $p<\infty$, which shall attract most of our attention. 
A finite measure $\mu$ on $M$ is said to have \emph{finite $p$-th moment} if for some
(hence all) $x_0\in M$ the following holds:
\[\int_{\mathbb{R}^d} d(x_0,x)^p \mu(\intd x) < +\infty.\]
In particular, any compactly supported finite measure has finite $p$-th moment for all $p$.

Let $\mu_0,\mu_1$ be two finite measures having finite $p$-th moment and the same mass.
A \emph{transport plan} between $\mu_0$ and $\mu_1$ is a measure $\Pi$ on $M\times M$
that has $\mu_0$ and $\mu_1$ as marginals, that is :
$\Pi(A\times M)=\mu_0(A)$ and $\Pi(M\times A)=\mu_1(A)$ for all Borelian set $A$.
One shall think of $\Pi$ has a assignement of mass: $\Pi(A\times B)$ represents
the mass sent from $A$ to $B$.

The $L^p$ \emph{cost} of a transport plan is defined as
\[c_p(\Pi) = \int_{M\times M} d(x,y)^p \Pi(\intd x\intd y).\]
One defines the $L^p$ \emph{Wasserstein distance} by
\[W_p(\mu_0,\mu_1)=\inf_\Pi c_p(\Pi)^{1/p}\]
where the infimum is on all tranport plan between $\mu_0$ and $\mu_1$. 
One can show that there is always a tranport plan that achieves this infimum, and
that $W_p$ defines a distance on the set of measures with finite $p$-th moment and given mass.

Moreover, if $M$ is compact $W_p$ metrizes the weak topology. 
If $M$ is non-compact, it defines a finer topology.

Most of the time, one restricts itself to probability measures. Here, we shall use
extensively mass transportation between submeasures of the main measures 
under study, so that
we need to consider measures of arbitrary mass. Given positive measures $\mu$ and $\nu$,
we write that $\mu\leqslant \nu$ if $\mu(A)\leqslant\nu(A)$ for all borelian set $A$, which means
that $\nu-\mu$ is also a positive measure.

It is important to
notice that $c_p(\Pi)$ is homogeneous of degree $1$ in the total mass and of degree
$p$ on distances, so that in the case $M=\mathbb{R}^d$ if $\varphi$ is a similitude of ratio $r$,
we have $W_p(m\, \varphi_\#\mu_0, m\, \varphi_\#\mu_1)=m^{1/p}\, r\, W_p(\mu_0,\mu_1)$.

The case $p=\infty$ is obtained as a limit of the finite case, see \cite{Champion-De_Pascale-Petri}.
Let $\mu_0$ and $\mu_1$ be compactly supported measures of the same mass and
let $\Pi$ be a transport plan between $\mu_0$ and $\mu_1$. The $L^\infty$ \emph{length}
of $\Pi$ is defined as
\[\ell_\infty(\Pi) = \sup\{d(x,y)\,|\, x,y\in\supp\Pi\}\]
that is, the maximal distance moved by some infinitesimal amount of mass when applying $\Pi$.
The $L^\infty$ distance between $\mu_0$ and $\mu_1$ then is
\[W_\infty(\mu_0,\mu_1)=\inf_\Pi \ell_\infty(\Pi)\]
where the infimum is on all transport plan from $\mu_0$ to $\mu_1$.
In a sense, the $L^\infty$ distance is a generalisation to measures of the Hausdorff metric
on compact sets. We shall use $\ell_\infty$, but not $d_\infty$. The problem of minimizing
$W_\infty(\mu,\Delta_N)$ is a matter of covering $\supp\mu$ (independently of $\mu$ itself),
a problem with quite a different taste than our.

\section{Preparatory results}

The following lemmas are useful tools we shall need; the first two at least
cannot pretend to any kind of originality by themselves.

\begin{lemm}[monotony]
Let $\mu$ and $\nu$ be finite measures of equal mass and $\tilde\mu\leqslant\mu$.
Then there is a measure $\tilde\nu\leqslant\nu$ 
(in particular, $\supp\tilde\nu\subset\supp\nu$) such that
\[W_p(\tilde\mu,\tilde\nu)\leqslant W_p(\mu,\nu).\]
\end{lemm}

\begin{proof}
Let $\Pi$ be an optimal transportation plan from $\mu$ to $\nu$.
We construct a low-cost transportation plan from $\tilde\mu$ to $\tilde\nu$
by disintegrating $\Pi$.

There is family of finite measures $(\eta_x)_{x\in M}$ such that
$\Pi = \int \eta_x \mu(\intd x)$, that is
\[\Pi(A\times B) = \int_{A} \eta_x(B) \mu(\intd x)\]
for all Borelian $A$ and $B$.
Define
\[\tilde\Pi(A\times B) = \int_A \eta_x(B) \tilde\mu(\intd x)\]
and let $\tilde\nu$ be the second factor projection of $\tilde\Pi$.
Since $\tilde\Pi\leqslant \Pi$, we have $\tilde\nu\leqslant \nu$
and $c_p(\tilde\Pi)\leqslant c_p(\Pi)$; moreover $\tilde\Pi$ is a transport plan 
from $\tilde\mu$ to $\tilde\nu$ by definition of $\tilde\nu$.
\end{proof}

\begin{lemm}[summing]
Let $(\mu,\nu)$ and $(\tilde\mu,\tilde\nu)$ be finite measures with pairwise
equal masses. Then
\[W_p^p(\mu+\tilde\mu,\nu+\tilde\nu)\leqslant W_p^p(\mu,\nu)+W_p^p(\tilde\mu,\tilde\nu).\]
\end{lemm}

\begin{proof}
Let $\Pi$ and $\tilde\Pi$ be optimal transport plans between respectively $\mu$ and $\nu$, $\tilde\mu$ and $\tilde\nu$.
Then $\Pi+\tilde\Pi$ is a transport plan between $\mu+\tilde\mu$ and $\nu+\tilde\nu$ whose cost
is $c_p(\Pi+\tilde\Pi)=c_p(\Pi)+c_p(\tilde\Pi)$.
\end{proof}

This very simple results have a particularly important consequence concerning our question.

\begin{lemm}[$L^1$ stability]
Let $\mu$ and $\tilde\mu$ be finite compactly supported measures on $M$,
$\varepsilon\in(0,1)$ and $(\mu_N)$ be any sequence of $N$-supported measures.

There is a sequence of $N$-supported measures $\tilde\mu_N$ such that  there are at most
$\varepsilon N$ points in $\supp\tilde\mu_N\setminus\supp\mu_N$ and
\[W_p^p(\tilde\mu,\tilde\mu_N)\leqslant W_p^p(\mu,\mu_{N_1}) 
+O\left(\frac{|\mu-\tilde\mu|_{TV}}{(\varepsilon N)^{p/d}}\right)\]
where $N_1$ is equivalent to (and at least) $(1-\varepsilon) N$, $|\cdot|_{TV}$ is the 
total variation norm and the constant in the $O$ depends only on the geometry of
a domain where both $\mu$ and $\tilde\mu$ are concentrated.

In particular we get \[W_p^p(\mu,\Delta_N)\leqslant W_p^p(\tilde\mu,\Delta_{N_1}) 
+O\left(\frac{|\mu-\tilde\mu|_{TV}}{(\varepsilon N)^{p/d}}\right).\]
\end{lemm}
The name of this result has been chosen to emphasize that the total variation distance
between two absolutely continuous measures is half the $L^1$ distance between their densities.

\begin{proof}
We can write $\tilde\mu=\mu' + \nu$ where $\mu'\leqslant\mu$ and $\nu$ is a positive
measure of total mass at most 
$|\mu-\tilde\mu|_{VT}$. If $D$ is a domain supporting  both $\mu$ and $\tilde\mu$,
it is a classical fact that there is a constant $C$ (depending only on $D$) such 
that for all integer $K$,
there are points  $x_1,\ldots,x_{K^d}\in D$ such that each point
of $D$ is at distance at most $C/K$ from one of the $x_i$. For example if
$D$ is a Euclidean cube of side length $L$, by dividing it regularly one
can achieve $C=L\sqrt{d}/2$.

Take $K=\lfloor(\varepsilon N)^{1/d}\rfloor$; then
by sending each point of $D$ to a closest $x_i$, one constructs a transport
plan between $\nu$ and a $K^d$-supported measure $\nu_N$ whose cost is at most 
$|\mu-\tilde\mu|_{VT} (C/K)^p$.

Let $N_1=N-K^d$. The monotony lemma gives a measure $\mu'_N\leqslant\mu_{N_1}$ (in particular, 
$\mu'_N$ is $N_1$-supported) such that
\[W_p(\mu',\mu'_N)\leqslant W_p(\mu,\mu_{N_1}).\]

The summing lemma now shows that
\[W_p^p(\tilde\mu,\mu'_N+\nu_N)\leqslant W_p^p(\mu,\mu_{N_1})+
  \frac{C^p |\mu-\tilde\mu|_{VT}}{K^p}.\]
\end{proof}
Note that the presence of the $|\mu-\tilde\mu|_{TV}$ factor will be crucial in the sequel,
but would not be present in the limit case $p=\infty$, which is therefore very different.

\begin{lemm}[metric stability]
Assume $D$ is a compact domain of $M$, endowed with
two different Riemannian metrics $g$ and $g'$ (defined on a open neighborhood of $D$).
Denote by $|g'-g|$ the minimum number $r$ such that 
\[e^{-2r} g_x(v,v)\leqslant g'_x(v)\leqslant e^{2r} g_x(v,v)\]
for all $x\in D$ and all $v\in T_x M$.

Then, denoting by $W_p$ the Wasserstein metric computed using the distance
$d$ induced by $g$, and by $W'_p$ the one obtained from the distance $d'$
induced by $g'$, one has for all measures $\mu,\nu$ supported on $D$ and of equal mass:
\[e^{-|g'-g|} W_p(\mu,\nu) \leqslant W'_p(\mu,\nu) \leqslant e^{|g'-g|} W_p(\mu,\nu).\]
\end{lemm}

\begin{proof}
For all $x,y\in D$ one has $d'(x,y)\leqslant e^r d(x,y)$
by computing the $g'$-length of a $g$-minimizing (or almost minimizing to avoid
regularity issues on the boundary) curve connecting $x$ to $y$. The same reasonning 
applies to transport plans:
if $\Pi$ is optimal from $\mu$ to $\nu$ according to $d$,
then the $d'$ cost of $\Pi$ is at most $e^{pr}$ times the $d$-cost of $\Pi$,
so that $W'_p(\mu,\nu) \leqslant e^r W_p(\mu,\nu)$. The other inequality follows
by symmetry.
\end{proof}

Let us end with a result showing that no mass is moved very far away by an optimal
transport plan to a $N$-supported measure if $N$ is large enough.
\begin{lemm}[localization]\label{lem:length}
Let $\mu$ be a compactly supported finite measure.
If $\mu_N$ is a closest $N$-supported measure to $\mu$ in $L^p$
Wasserstein distance
and $\Pi_N$ is a $L^p$ optimal transport plan between $\mu$ and $\mu_N$,
then when $N$ goes to $\infty$,
\[\ell_\infty(\Pi_N)\to 0.\]
\end{lemm}

\begin{proof}
Assume on the contrary that there are sequences $N_k\to\infty$,
$x_k\in\supp\mu$ and a number $\varepsilon>0$ such that
$\Pi_{N_k}$ moves $x_k$ by a distance at least $\varepsilon$.
There is a covering of $\supp\mu$ by a finite number of balls of radius
$\varepsilon/3$. Up to extracting a subsequence, we can assume that all
$x_k$ lie in one of this balls, denoted by $B$. Since $B$ is a neighborhood
of $x_k$ and $x_k\in\supp\mu$, we have $\mu(B)>0$. Since $\Pi_{N_k}$ is
optimal, it moves $x_k$ to a closest point in $\supp\mu_{N_k}$, which must be at distance
at least $\varepsilon$ from $x_k$. Therefore, every point in $B$ is at distance
at least $\varepsilon/3$ from $\supp\mu_{N_k}$, so that
$c_p(\Pi_{N_k})\geqslant \mu(B) (\varepsilon/3)^p >0$, in contradiction
with $W_p(\mu,\Delta_N)\to 0$.
\end{proof}

\section{Approximation rate and dimension}\label{sec:dimension}

Theorem \ref{theo:dimension} is the union of the two following propositions.
Note that the estimates given do not depend much on $p$, so that in fact
Theorem \ref{theo:dimension} stays true when $p=\infty$.

\begin{prop}\label{prop:majo}
If $\mu$ is a compactly supported probability measure on $M$ and if for some $C>0$
and for all $r\leqslant\diam(\supp\mu)$, one has 
\[C^{-1} r^s\leqslant \mu(B(x,r))\]
then  for all $N$
\[W_p(\mu,\Delta_N)\leqslant \frac{5C^{1/s}}{N^{1/s}}.\]
\end{prop}

\begin{proof}
The $5\delta$ covering proposition above implies that given any $\delta>0$, there is a subset
$\mathscr{G}$ of
$\supp\mu$ such that 
\begin{itemize}
\item $\supp\mu \subset \bigcup_{x\in \mathscr{G}} B(x,5\delta)$,
\item $B(x,\delta)\cap B(x',\delta)=\varnothing$ whenever $x\neq x' \in\mathscr{G}$.
\end{itemize}
In particular, as soon as $\delta<\diam(\supp\mu)$ one has
\[1\geqslant \sum_{x\in\mathscr{G}} \mu(B(x,\delta)) \geqslant |\mathscr{G}| C^{-1} \delta^s\]
so that $\mathscr{G}$ is finite, with $ |\mathscr{G}|\leqslant C\delta^{-s}$.

Let $\tilde\mu$ be a measure supported on $\mathscr{G}$, that minimizes the $L^p$ distance
to $\mu$ among those. A way to construct $\tilde\mu$ is to assign to a point $x\in\mathscr{G}$
a mass equal to the $\mu$-measure of its Voronoi cell, that is of the set of points
nearest to $x$ than to any other points in $\mathscr{G}$. The mass at a point at equal
distance from several elements of $\mathscr{G}$ can be split indifferently between those.
The previous discussion also gives a transport plan from $\mu$ to $\tilde\mu$,
where each bit of mass moves a distance at most $5\delta$, so that
$W_p(\mu,\tilde\mu)\leqslant 5\delta$ (whatever $p$).

Now, let $N$ be a positive integer and choose $\delta = (C/N)^{1/s}$. The family $\mathscr{G}$
obtained from that $\delta$ has less than $N$ elements, so that
$W_p(\mu,\Delta_N)\leqslant 5(C/N)^{1/s}$.
\end{proof}

\begin{prop}
If $\mu$ is a probability measure on $M$ and if for some $C>0$
and for all $r$, one has 
\[ \mu(B(x,r))\leqslant C r^s \]
then for all $N$,
\[ \left(\frac{s}{s+p}\right)^{1/p} C^{-1/s}\ \frac1{N^{1/s}}\leqslant W_p(\mu,\Delta_N).\]
\end{prop}

\begin{proof}
Consider a measure $\mu_N\in\Delta_N$ that minimizes the distance
to $\mu$. For all $\delta>0$, the union of the balls centered at $\supp\mu_N$ and of radius
$\delta$ has $\mu$-measure at most $NC\delta^s$. In any transport plan from $\mu$ to
$\mu_N$, a quantity of mass at least $1-NC\delta^s$ travels a distance at least $\delta$,
so that in the best case the quantity of mass traveling a distance between $\delta<(NC)^{-1/s}$
and $\delta+\intd\delta$ is $NCs\delta^{s-1}\intd \delta$. It follows that
\[W_p(\mu,\mu_N)^p\geqslant \int_0^{(NC)^{-1/s}} sNC\delta^{s-1}\delta^p \intd \delta\]
so that $W_p(\mu,\mu_N)\geqslant (s/(s+p))^{1/p}(NC)^{-1/s}$.
\end{proof}

In fact, Theorem \ref{theo:dimension} applies to more general measures, for example
combination of Ahlfors regular ones, thanks to the following.
\begin{lemm}
If $\mu=a_1\mu^1+a_2\mu^2$ where $a_i>0$ and $\mu^i$ are probability measures
such that 
$W_p(\mu^2,\Delta_N)\lesssim W_p(\mu^1,\Delta_N)$
and $W_p(\mu^1,\Delta_N)\lesssim W_p(\mu^1,\Delta_{2N})$
 then
$W_p(\mu,\Delta_N)\approx W_p(\mu^1,\Delta_N)$.
\end{lemm}

\begin{proof}
By the monotony lemma, $W_p(a_1\mu^1,\Delta_N)\leqslant W_p(\mu,\Delta_N)$ so that
$W_p(\mu^1,\Delta_N)\leqslant a_1^{-1/p}W_p(\mu,\Delta_N)$.

The summing lemma gives
\[W_p(\mu,\Delta_{2N})\leqslant\big(W_p(a_1\mu^1,\Delta_N)^p+W_p(a_2\mu^2,\Delta_N)^p\big)^{1/p}\]
so that 
\[W_p(\mu,\Delta_{2N})\lesssim W_p(a_1\mu^1,\Delta_N)\lesssim W_p(\mu^1,\Delta_{2N}).\]
Since $W_p(\mu,\Delta_{2N+1})\leqslant W_p(\mu,\Delta_{2N})$ we also get
\[W_p(\mu,\Delta_{2N+1})\lesssim W_p(\mu^1,\Delta_{2N}) \lesssim W_p(\mu^1,\Delta_{4N})\leqslant W_p(\mu^1,\Delta_{2N+1})\]
\end{proof}

The following is now an easy consequence of this lemma.
\begin{coro}
Assume that $\mu=\sum_{i=1}^k a_i \mu^i$ where
$a_i>0$ and $\mu^i$ are probability measures that
are compactly supported and Ahlfors regular of dimension $s_i>0$.
Let $s=\max_i(s_i)$. Then
\[W_p(\mu,\Delta_N)\approx \frac1{N^{1/s}}.\]
\end{coro}

\section{Absolutely continuous measures}\label{sec:cont}

In this section we prove Theorem \ref{theo:cont}.
To prove the Euclidean case, the idea is to approximate (in the $L^1$ sense)
a measure with density by a combination of uniform measures in squares. Then
a measure on a manifold can be decomposed as a combination of measures supported in
charts, and by metric stability the problem reduces to the Euclidean case.

The following key lemma shall be used several times to extend the class
of measures for which we have precise approximation estimates.
\begin{lemm}[Combination]\label{lemm:combination}
Let $\mu$ be an absolutely continuous measure on $M$. 
Let $D_i$ ($1\leqslant i\leqslant I$) be domains of $M$ whose interiors
do not overlap, and assume we can decompose $\mu=\sum_{i=1}^I \mu^i$
where $\mu^i$ is non-zero and supported on $D_i$.
Assume moreover that there are
numbers $(\alpha_1,\ldots,\alpha_I)=\alpha$ such that
\[W_p(\mu^i,\Delta_N)\sim \frac{\alpha_i}{N^{1/d}}.\]

Let $\mu_N\in\Delta_N$ be a sequence minimizing the $W_p$ distance to $\mu$
and define $N_i$ (with implicit dependency on $N$)
as the number of points of $\supp\mu_N$ that lie on $D_i$,
the points lying on a common boundary of two or more domains being attributed
arbitrarily to one of them.

If the vector $(N_i/N)_i$ has a cluster point $x=(x_i)$, then
$x$ minimizes 
\[F(\alpha;x)=\left(\sum_i \frac{\alpha_i^p}{x_i^{p/d}}\right)^{1/p}\] 
and if $(N_i/N)_i\to x$ when $N\to\infty$, then
\[W_p(\mu,\mu_N)\sim \frac{F(\alpha;x)}{N^{1/d}}.\]
\end{lemm}
Note that the assumption that none of the $\mu^i$ vanish is obviously
unecessary (but convenient). If some of the $\mu^i$ vanish, one only has to
dismiss them. 

\begin{proof}
For simplicity we denote $c_p(N)=W_p^p(\mu,\Delta_N)$.
Let $\varepsilon$ be any positive, small enough number.

We can find a $\delta>0$ and domains $D'_i\subset D_i$
such that: each point of $D'_i$ is at distance at least $\delta$
from the complement of $D_i$; if $\mu'^i$ is the restriction
of $\mu^i$ to $D'_i$, $|\mu'^i-\mu^i|_{VT}\leqslant \varepsilon^{1+p/d}$.

Assume $x$ is the limit of $(N_i/N)_i$ when $N\to \infty$.
Let us first prove that none of the $x_i$ vanishes.
Assume the contrary for some index $i$: then $N_i=o(N)$.
For each $N$ choose an optimal transport plan $\Pi_N$ from $\mu$
to $\mu_N$. 
Let $\nu_{N}\leqslant \mu^i$ be the part of $\mu^i$
that is sent by $\Pi_N$ to the $N_i$ points of $\supp\mu_N$ that lie
in $D_i$, constructed as in the summing lemma,
and let $m_N=\mu^i(D'_i)-\nu_N(D'_i)$ be the mass
that moves from $D'_i$ to the exterior of $D_i$ under $\Pi_N$.
Then the cost of $\Pi_N$ is bounded from below by 
$m_N \delta^p+W_p^p(\nu_N,\Delta_{N_i})$. Since it goes to zero,
we have $m_N\to0$ and up to extracting a subsequence $\nu_N\to\nu$
where $\mu'^i\leqslant\nu\leqslant \mu^i$. The cost of $\Pi_N$
is therefore bounded from below by all number less than
$W_p^p(\nu_N,\Delta_{N_i})\lesssim N_i^{-1/d} \ll N^{-1/d}$,
a contradiction.

Now, let $\varepsilon$ be any positive, small enough number.
By considering optimal transport plans between $\mu^i$
and optimal $N_i$-supported measures of $D_i$, we get that
\begin{eqnarray*}
c_p(N) &\leqslant& \sum_i W_p^p(\mu^i,\Delta_{N_i})\\
       &\leqslant& \sum_i \frac{(\alpha_i+\varepsilon)^p}{N_i^{p/d}}
\end{eqnarray*}
when all $N_i$ are large enough, which happens if $N$ itself is large enough given
that $x_i\neq0$.

For $N$ large enough, the localization lemma ensures that no mass is
moved more than $\delta$ by an optimal transport plan between
$\mu$ and $\mu_N$. This implies that the cost $c_p(N)$ is bounded
below by $\sum_i W_p^p(\mu'^i,\Delta_{N_i})$. By $L^1$-stability
this gives the bound
\[c_p(N)\geqslant \frac{\alpha_i^p(1-\varepsilon)^{p/d}+O(\varepsilon)}{(x_iN)^{p/d}}.\]

The two inequalities above give us
\[c_p(N) N^{p/d} \to \sum_i \frac{\alpha_i^p}{x_i^{p/d}}=F^p(\alpha;x).\]
Now, if $x$ is a mere cluster point of $(N_i/N)$, this
still holds up to a subsequence. If $x$ did not minimize
$F(\alpha;x)$, then by taking best approximations of $\mu^i$
supported on $x'_iN$ points where $x'$ is a minimizer, we would get by the same computation 
a sequence $\mu'_N$ with better asymptotic behavior than $\mu_N$
(note that we used the optimality of $\mu_N$ only to bound from above
each $W_p^p(\mu^i,\Delta_{N_i})$).
\end{proof}

The study of the functional $F$ is straightforward.
\begin{lemm}
Fix a positive vector $\alpha=(\alpha_1,\alpha_2,\ldots,\alpha_I)$ and consider the simplex
$X=\{(x_i)_i\,|\, \sum_i x_i =1,\, x_i> 0\}$. The function
$F(\alpha;\cdot)$  has a unique minimizer $x^0=(x_i^0)$ in $X$, 
which is proportionnal to $(\alpha_i^{\frac {dp}{d+p}})_i$, with
\[F(\alpha ; x^0)=\left(\sum_i \alpha_i^{\frac {dp}{d+p}} \right)^{\frac{d+p}{dp}}
  =:|\alpha|_{\frac{dp}{d+p}}.\]

As a consequence, in the combination lemma the vector $(N_i/N)$ 
must converge to $x^0$.
\end{lemm}

\begin{proof}
First $F(\alpha;\cdot)$ is continuous and goes to $\infty$ on the boundary of $X$,
so that is must have
a minimizer. Any minimizer must be a critical point of $F^p$ and therefore satisfy
\[\sum_i \alpha_i^p x_i^{-p/d-1} \eta_i=0\]
for all vector $(\eta_i)$ such that $\sum_i \eta_i=0$. This holds only when
 $\alpha_i^p x_i^{-p/d-1}$ is a constant and we get the uniqueness of $x^0$ and its expression:
\[x^0_i = \frac{\alpha_i^{\frac{dp}{d+p}}}{\sum_j\alpha_j^{\frac{dp}{d+p}}}.\]
The value of $F(\alpha; x^0)$ follows.

In the combination lemma, we now by compacity that $(N_i/N)$ must have cluster points,
all of which must minimize $F(\alpha;\cdot)$. Since there is only one minimizer,
$(N_i/N)$ has only one cluster point and must converge to $x^0$.
\end{proof}

We are now ready to tackle more and more cases in Theorem \ref{theo:cont}. 
As a starting point, we consider the uniform measure $\square^d$ on the unit 
cube of $\mathbb{R}^d$ (endowed with the canonical metric).
\begin{prop}\label{prop:cube}
There is a number $\theta(d,p)>0$ such that
\[W_p(\square^d,\Delta_N)\sim \frac{\theta(d,p)}{N^{1/d}}.\]
\end{prop}

The proof is obviously not new, since it is the same argument that shows
that an optimal packing (or covering) of the Euclidean space must have a well-defined
density (its upper and lower densities are equal).
\begin{proof}
Let $c(N)=W_p^p(\square^d,\Delta_N)$.
We already know that $c(N)\approx N^{-p/d}$, so let
$A=\liminf N^{p/d} c(N)$ and consider any $\varepsilon>0$.

Let $N_1$ be an integer such that $c(N_1)\leqslant (A+\varepsilon) N_1^{-p/d}$
and let $\mu_1\in \Delta_{N_1}$ be nearest to $\mu$.
For any integer $\ell$, we can write $\ell=k^d+q$ where $k=\lfloor\ell^{1/d}\rfloor$
and $q$ is an integer; then $q = O(\ell^{1-1/d})=o(\ell)$ where the $o$ depends only on $d$.

\begin{figure}[tp]\begin{center}
\includegraphics{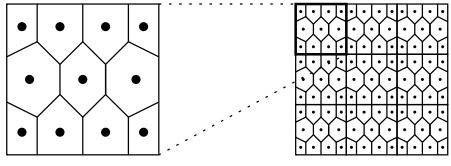}
\caption{An optimal $N$-supported measure can be used to construct a good
$k^dN$-supported measure for all $k$.}\label{fig:cube}
\end{center}\end{figure}

Divide the cube into $k^d$ cubes of side length $1/k$, and consider the element
$\mu_k$ of $\Delta_{k^dN_1}$ obtained by duplicating $\mu_1$ in each of the cubes,
with scaling factor $k^{-1}$ and mass factor $k^{-d}$ (see figure \ref{fig:cube}).
The obvious transport plan obtained
in the same way from the optimal one between $\square^d$ and $\mu_1$ has total 
cost $k^{-p}c(N_1)$, so that
\[c(\ell N_1)\leqslant k^{-p} \frac{A+\varepsilon}{N_1^{p/d}}=
  \left(\frac\ell {k^d}\right)^{p/d} \frac{A+\varepsilon}{(\ell N_1)^{p/d}}.\]
But since $k^d\sim\ell$, for $\ell$ large enough we get
\[c(\ell N_1)\leqslant\frac{A+2\varepsilon}{(\ell N_1)^{p/d}}.\]
Now $N\sim \lfloor N/N_1\rfloor N_1$ so that for $N$ large enough
$c(N)\leqslant(A+3\varepsilon) N^{-p/d}$.
This proves that $\limsup N^{p/d} c(N)\leqslant A+3\varepsilon$ for all
$\varepsilon>0$.
\end{proof}
Note that we used the self-similarity of the cube at many different scales;
the result does not hold with more general self-similar (fractal) measures,
see Section \ref{sec:examples}.

Now, the combination lemma enables us to extend the validity domain of
Equation \eqref{eq:main}.
\begin{lemm}
Let $\mu=\rho\lambda$ where $\lambda$ is the Lebesgue measure on
$\mathbb{R}^d$, $\rho$ is a $L^1$ non-negative function supported on
a union of cubes $C_i$ with non-overlapping interiors, side length $\delta$,
and assume $\rho$ is constant on each cube, with value $\rho_i$.
Then Equation \eqref{eq:main} holds.
\end{lemm}

\begin{proof}
Let $\mu^i$ the restriction of $\mu$ to $C_i$, removing any cube where $\rho$
vanishes identically. Then from Proposition \ref{prop:cube} we get
$W_p(\mu^i,\Delta_N)\sim \alpha_i N^{-1/d}$ where 
\[\alpha_i=\theta(d,p)(\rho_i\delta^d)^{1/p} \delta= \theta(d,p) \rho_i^{1/p} \delta^{\frac{d+p}p}\] 
due to the homogeneity of $W_p$:
$\mu_i$ is obtained from $\square^d$ by multiplication by $\rho_i\delta^d$
and dilation of a factor $\delta$. By the combination lemma, we get
$W_p(\mu,\Delta_N)\sim \min F(\alpha;\cdot) N^{-1/d}$
where 
\begin{eqnarray*}
\min F(\alpha,\cdot)&=&\theta(d,p)\left|\sum_i \rho_i^{\frac d{d+p}} 
                       \delta^d\right|^{\frac{d+p}{dp}}\\
  &=& \theta(d,p) |\rho|_{\frac d{d+p}}^{1/p}.
\end{eqnarray*}
\end{proof}

\begin{lemm}
Equation \eqref{eq:main} holds whenever $\mu$ is an absolutely
continuous measure defined on a compact domain of $\mathbb{R}^d$.
\end{lemm}

\begin{proof}
For simplicity, we denote $\beta=d/(d+p)$.
Let $C$ be a cube containing the support of $\mu$.
Choose some $\varepsilon>0$. Let $\tilde\mu=\tilde\rho\lambda$ be a measure such that
$\tilde\rho$ is constant on each cube of a regular subdivision of $C$, is zero outside $C$,
satisfies $|\rho-\tilde\rho|_1\leqslant 2\varepsilon^{1+p/d}$ and
such that $|\rho-\tilde\rho|_\beta\leqslant \varepsilon|\rho|_\beta$.

The stability lemma shows that
\[W_p^p(\mu,\Delta_N) \leqslant W_p^p(\tilde\mu,\Delta_{(1-\varepsilon)N}) 
       +O\left(\frac{|\rho-\tilde\rho|_1}{2(\varepsilon N)^{p/d}}\right)\]
so that, using the hypotheses on $\tilde\rho$ and the previous lemma,
\[W_p^p(\mu,\Delta_N) \leqslant 
\frac{(\theta(d,p)+\varepsilon)^p |\rho|_\beta (1+\varepsilon)(1-\varepsilon)^{-p/d}
+O(\varepsilon)}{N^{p/d}}\]
for $N$ large enough.

Symmetrically, we get (again for $N$ large enough)
\begin{eqnarray*}
W_p^p(\mu,\Delta_N) &\geqslant& W_p^p(\tilde\mu,\Delta_{N/(1-\varepsilon)}) 
      -O\left(\frac{|\rho-\tilde\rho|_1}{2(\varepsilon N)^{p/d}}\right) \\
     &\geqslant& 
\frac{(\theta(d,p)-\varepsilon)^p|\rho|_\beta(1-\varepsilon)^{1+p/d}-O(\varepsilon)}{N^{p/d}}.
\end{eqnarray*}
Letting $\varepsilon\to0$, the claimed equivalent follows.
\end{proof}

\begin{lemm}
Equation \eqref{eq:main} holds whenever $\mu$ is an absolutely
continuous measure defined on a compact domain of $\mathbb{R}^d$, 
endowed with any Riemannian metric.
\end{lemm}

\begin{proof}
Denote by $g$ the Riemannian metric, and let $C$ be a Euclidean cube 
containing the support of $\mu$. Let $\varepsilon$ be any positive number,
and choose a
regular subdivision of $C$ into cubes $C_i$ of center $p_i$ such that for all $i$,
the restriction $g_i$ of $g$ to $C_i$ is almost constant:
$|g(p)-g(p_i)|\leqslant \varepsilon/2$ for all $p\in C_i$. Denote
by $\tilde g$ the piecewise constant metric with value $g(p_i)$ on $C_i$.
Note that even if $\tilde g$ is not continuous, at each discontinuity
point $x$ the possible choices for the metric are within a factor
$e^{2\varepsilon}$ one from another, and one defines that $\tilde g(x)(v,v)$
is the least of the $g(p_i)(v,v)$ over all $i$ such that $x\in C_i$. In this
way, $\tilde g$ defines a distance function close to the distance induced
by $g$ and the metric stability lemma holds with the same proof.

If one prefers not using discontinuous metrics, then it is also possible to
consider slightly smaller cubes $C'_i\subset C_i$, endow $C'_i$ with a constant
metric, and interpolate the metric between the various cubes. Then one uses
the $L^1$ stability in addition to the metric stability in the sequel.

Denote by $\rho$ the density of $\mu$ with respect to the volume form defined by
$g$, by $\mu^i$ the restriction of $\mu$ to $C_i$ and by $\rho_i$ the density of
$\mu^i$.
A domain of $\mathbb{R}^d$ endowed with a constant metric is isometric
to a domain of $\mathbb{R}^d$ with the Euclidean metric so that we can apply
the preceding lemma to each $\mu^i$: denoting by $W'_p$ the Wasserstein distance
computed from the metric $\tilde g$,
\[ W'_p(\mu^i,\Delta_N)\sim\frac{\delta(d,p) |\rho_i|_{\frac d{d+p}}^{1/p}}{N^{1/d}}.\]
The combination lemma then ensures that
$W'_p(\mu,\Delta_N)\sim \min F(\alpha;\cdot) N^{-1/d}$
where
\begin{eqnarray*}
\min F(\alpha,\cdot)&=&\theta(d,p)\left(\sum_i \int_{C_i} \rho_i^{\frac d{d+p}} 
                       \right)^{\frac{d+p}{dp}}\\
  &=& \theta(d,p) |\rho|_{\frac d{d+p}}^{1/p}.
\end{eqnarray*}
The metric stability lemma gives 
\[e^{-\varepsilon} W'_p(\mu,\Delta_N) \leqslant W_p(\mu,\Delta_N)
  \leqslant e^\varepsilon W'_p(\mu,\Delta_N)\]
and we only have left to let $\varepsilon\to 0$.
\end{proof}

We can finally end the proof of the main theorem.
\begin{proof}[Proof of Theorem \ref{theo:cont}]
Here $\mu$ is an absolutely
continuous measure defined on a compact domain $D$ of $M$.
Divide the domain into a finite number of subdomains $D_i$, each of which
is contained in a chart. Using this chart, each $D_i$ is identified with
a domain of $\mathbb{R}^d$ (endowed with the pulled-back metric of $M$).
By combination, the previous lemma shows that Equation \eqref{eq:main}
holds.

Let us now give the asymptotic distribution of the support of any
distance minimizing $\mu_N$. Let $A$ be any domain in $M$.
Let $x$ be the limit of the proportion of $\supp\mu_N$ that lies inside $A$ ($x$ exists up
to extracting a subsequence). Since the domains generates the Borel $\sigma$-algebra, we
only have to prove that $x=\int_A \rho^\beta / \int_M \rho^\beta$. But this
follows from the combination lemma applied to the restriction of $\mu$
to $A$ and to its complement.
\end{proof}

\section{The dyadic Cantor measure}\label{sec:examples}

In this section we study the approximation problem for the dyadic Cantor measure
$\kappa$ to prove Theorem \ref{theo:example}.

Let $S^0, S^1$ be the dilations of ratio $1/3$ and fixed point $0,1$.
The map defined by
\[\mathscr{S}:\mu\mapsto 1/2\, S^0_\#\mu + 1/2\, S^1_\#\mu\]
is $1/3$-Lipschitz on the complete metric space of probability measures
having finite $p$-th moment endowed with the $L^p$ Wasserstein metric.
It has therefore a unique fixed point, called the dyadic Cantor measure
and denoted by $\kappa$. It can be considered as the ``uniform'' measure on
the usual Cantor set.

By convexity of the cost function and symmetry, $c_1:=W_p(\kappa,\Delta_1)$
is realized by the Dirac measure at $1/2$. Using the contractivity
of $\mathscr{S}$, we see at once that $W_p(\kappa,\Delta_{2^k})\leqslant 3^{-k} c_1$.
Denote by $s=\log 2/\log 3$ the dimension of $\kappa$. We have
\[ W_p(\kappa,\Delta_{2^k})(2^k)^{1/s}\leqslant c_1\]
for all integer $k$.

To study the case when the number of points is not a power of $2$, and to get
lower bounds in all cases, we introduce a notation to code the regions of $\supp\kappa$.
Let $I^0=[0,1]$ and given a word $w=\epsilon_n\ldots\epsilon_1$ where
$\epsilon_i\in\{0,1\}$, define $I_w^n = S_{\epsilon_n} S_{\epsilon_{n-1}} \cdots S_{\epsilon_1} [0,1]$.
The \emph{soul} of such an interval is the open interval of one-third length with the same center.
The \emph{sons} of $I_w^n$ are the two intervals $I_{\epsilon w}^{n+1}$
where $\epsilon\in\{0,1\}$, and an interval is the \emph{father} of its sons.
The two sons of an interval are \emph{brothers}. Finally,
we say that $n$ is the \emph{generation} of the interval $I_w^n$.

Let $N$ be an integer, and $\mu_N\in\Delta_N$ be a measure closest to $\kappa$,
whose support is denoted by $\{x_1,\ldots,x_N\}$.
An interval $I_w^n$ is said to be \emph{terminal} if there is an $x_i$ in its soul.
A point in $I_w^n$ is always closer to the center of $I_w^n$ than to the center of its
father. This and the optimality of $\mu_N$ implies that a terminal interval contains only one $x_i$, at its center. 

Since the restriction of $\kappa$ to $I_w^n$ is a copy of $\kappa$
with mass $2^{-n}$ and size $3^{-n}$, it follows that 
\[W_p(\kappa,\mu_N)^p=W_p^p\sum_{I_w^n} 2^{-n}3^{-np}\]
where the sum is on terminal intervals. A simple convexity arguments shows that the terminal intervals
are of at most two (successive) generations.

\begin{figure}[tp]\begin{center}
\input{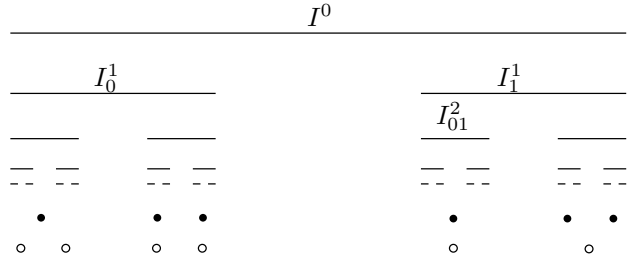}
\caption{The four first steps of the construction of the Cantor set; the Cantor measure is 
equally divided between the intervals of a given step. The bullets show the supports
of two optimal approximation of $\kappa$ by $6$-supported measures. We see that there is
no need for the support to be equally distributed between the intervals of the first 
generation.}\label{fig:cantor}
\end{center}\end{figure}

Consider the numbers $N_k=3\cdot2^k$. The terminal intervals of an optimal $\mu_{N_k}$
must be in generations $k+1$ (for $2^k$ of them) and $k+2$ (for $2^{k+1}$ of them).
Therefore
\[W_p(\kappa,\mu_{N_k})^p=c_1^p\left(3^{-(k+1)p}+3^{-(k+2)p} \right)/2\]
and finally
\[W_p(\kappa,\Delta_{N_k}) N_k^{1/s} = c_1 \left(\frac{1+3^{-p}}2\right)^{1/p} 3^{\frac{\log3}{\log2}-1}.\]
Note that the precise repartition of the support does not have any importance (see figure 
\ref{fig:cantor}).

To see that $W_p(\kappa,\Delta_N)N^{1/s}$ has no limit, it is now sufficient to estimate
the factor of $c_1$ in the right-hand side of the above formula. First we remark that
$\left(\frac{1+3^{-p}}2\right)^{1/p}$ is greater than
$1-(1-3^{-p})/(2p)$ which is increasing in $p$ and takes for $p=1$ the value
$2/3$. Finally, we compute $2/3\cdot3^{\frac{\log3}{\log2}-1}\simeq 1.27 >1$.

Note that the fundamental property of $\kappa$ we used is that the points in a
given $I_w^n$ are closest to its center than to that of its father. The same
method can therefore be used to study the approximation of sparser Cantor measure, or
to some higher-dimensionnal analogue like the one generated by four
contractions of ratio $1/4$ on the plane, centered at the four vertices of a square.

Moreover, one could study into more details the variations in the approximations
$W_p(\kappa,\Delta_N)$. As said before, here our point was only to show the limitations
to Theorem \ref{theo:cont}.

\section{Link with Centroidal Voronoi Tessellations}\label{sec:CVT}

Here we explain the link between our optimization problem and the centroidal Voronoi tessellations
(CVTs in short).
For a complete account on CVTs, the reader can consult
\cite{Du-Faber-Gunzburger} from where all definitions below
are taken. Since we use the concept of barycenter, we consider only the case $M=\mathbb{R}^d$
(with the Euclidean metric). As before, $\lambda$ denotes the Lebesgue measure.

\subsection{A quick presentation}

Consider a compact convex domain $\Omega$ in $\mathbb{R}^d$ and a density (positive, $L^1$)
function $\rho$ on $\Omega$.

Given a $N$-tuple $X=(x_1,\ldots,x_N)$ of so-called \emph{generating points}, one defines
the associated \emph{Voronoi Tessellation} as the collection of convex sets
\[V_i = \big\{x\in\Omega\,\big|\, |x-x_i|\leqslant |x-x_j|\mbox{ for all } j\in\llbracket 1,N \rrbracket \big\}\] 
and we denote it by $V(X)$. One says that $V_i$ is the \emph{Voronoi cell} of
$x_i$. It is a tiling of $\Omega$, in particular the cells cover $\Omega$
and have disjoint interiors.

Each $V_i$ has a center of mass, equivalently defined as 
\[g_i=\frac{\int_{V_i} x\rho(x) \intd x}{\int_{V_i} \rho(x)\intd x}\]
or as the minimizer of the energy functionnal 
\[\mathscr{E}_{V_i}(g)=\int_{V_i} |x-g|^2\rho(x)\intd x.\]

One says that $(V_i)_i$ is a \emph{centroidal Voronoi tessellation} or \emph{CVT},
if for all $i$, $g_i=x_i$. The existence of CVTs comes easily by considering the following
optimization problem: search for a $N$-tuple of points $X=(x_1,\ldots,x_N)$
and a tiling $V$ of $\Omega$ by $N$ sets $V_1,\ldots,V_N$ which together minimize
\[\mathscr{E}_V(X)=\sum_{i=1}^N \mathscr{E}_{V_i}(x_i)\]
A compacity argument shows that such a minimizer exists, so
let us explain why a minimizer must be a CVT together with its generating set.
First, each $x_i$ must be the center of mass of $V_i$, otherwise one could reduce
the total energy by moving $x_i$ to $g_i$ and changing nothing else. But also, $V_i$ should
be the Voronoi cell of $x_i$, otherwise there is a $j\neq i$ and a set of positive measure in $V_i$
whose points are closest to $x_j$ than to $x_i$. Transfering this set from $V_i$ to $V_j$
would reduce the total cost.

 We observe that
this optimization problem is exactly that of approximating the measure
$\rho\lambda$ in $L^2$ Wasserstein distance; more precisely,
finding the $N$-tuple $x$ that minimizes $\inf_V \mathscr{E}_V(X)$ is equivalent
to finding the support of an optimal $\mu_N\in\Delta_N$ closest to $\rho\lambda$,
and then the Voronoi tesselation generated by $X$ gives the mass of $\mu_N$
at each $x_i$ and the optimal transport from $\rho\lambda$ to $\mu_N$.

One says that a CVT is \emph{optimal} when its generating set is a global minimizer
of the energy functional
\[\mathscr{E}(X)=\mathscr{E}_{V(X)}(X).\]
Optimal CVTs are most important in applications, which include
for example mesh generation and image analysis (see \cite{Du-Faber-Gunzburger}).

\subsection{Equidistribution of Energy}

The \emph{principle of energy equidistribution} says that if $X$ generates an optimal
CVT, the energies $\mathscr{E}_{V_i}(x_i)$ of the generating points should be asymptotically
independent of $i$ when $N$ goes to $\infty$.

Our goal here is to deduce a mesoscopic version of this principle from Theorem \ref{theo:cont}.
A similar result holds for any exponent, so that we introduce the $L^p$ energy functionals
$\mathscr{E}^p_{V_i}(x_i) = \int_{V_i} |x-x_i|^p\rho(x)\intd x$,
$\mathscr{E}^p_V(X) = \sum_i \mathscr{E}^p_{V_i}(x_i)$ and
$\mathscr{E}^p(X) = \inf_V \mathscr{E}^p_V(x) = \mathscr{E}^p_{V(X)}(X)$.
In particular, an optimal $X$ for this last functional is the support
of an element of $\Delta_N$ minimizing the $L^p$ Wasserstein distance to $\rho\lambda$.
  
Note that for $p\neq2$ an $x$ minimizing $\mathscr{E}^p(x)$ need not generate
a CVT, since the minimizer of $\mathscr{E}^p_{V_i}$ is not always
the center of mass of $V_i$ (but it is unique as soon as $p>1$).

\begin{coro}\label{coro:equid}
Let $A$ be a cube of $\Omega$.
Let $X^N=\{x_1^N,\ldots,x_N^N\}$ be a sequence
of $N$-sets minimizing $\mathscr{E}^p$ for the density $\rho$, and denote
by $\bar{\mathscr{E}}^p_A(N)$ the average energy of the points of $X^N$ that lie in $A$.
Then
\[\bar{\mathscr{E}}^p_A(N) N^{\frac{d+p}d}\]
has a limit when $N\to\infty$, and this limit does not depend on $A$.
\end{coro}
The cube $A$ could be replaced by any domain,
but not by any open set. Since the union of the $X_N$ is countable, there are
indeed open sets of arbitrarily small measure containing all the points $(x_i^N)_{N,i}$.


\begin{proof}
Fix some $\varepsilon>0$ and let
$A'\subset A$ be the set of points that are at distance at least
$\varepsilon$ from $\Omega\setminus A$ and by $A''\supset A$ the set of points
at distance at most $\varepsilon$ from $A$.

First, the numbers $N',N''$ of points of $X^N$ in $A',A''$ satisfy
\[N'\sim N \frac{\int_{A'} \rho^{d/(d+p)}}{\int_\Omega \rho^{d/(d+p)}}\qquad 
N''\sim N \frac{\int_{A''} \rho^{d/(d+p)}}{\int_\Omega \rho^{d/(d+p)}}.\]

The localization lemma implies that the maximal distance by which mass is moved by the optimal
transport between $\rho\lambda$ and the optimal $X^N$-supported
measure tends to $0$, so that for $N$ large enough the energy of all points in $A$
is at least the minimal cost between $\rho_{|A'}\lambda$ and 
$\Delta_{N'}$ and at most
the minimal cost between $\rho_{|A''}\lambda$ and 
$\Delta_{N''}$.

Letting $\varepsilon\to 0$ we thus get
that the total energy of all points of $X^N$ lying in $A$ is equivalent
to
\[\theta(d,p)\frac{\left(\int_A \rho^{d/(d+p)}\right)^{(d+p)/d}}{\left(N\int_A \rho^{d/(d+p)}\right)^{p/d}}
= \theta(d,p) N^{-p/d}\int_A \rho^{d/(d+p)} \]

As a consequence we have
$\bar{\mathscr{E}}_A(N)\sim  \left(\theta(d,p) \int_\Omega\rho^{d/(d+p)}\right)  N^{-(d+p)/d}$.
\end{proof}

\bibliographystyle{smfplain}
\bibliography{biblio}

\providecommand{\bysame}{\leavevmode ---\ }
\providecommand{\og}{``}
\providecommand{\fg}{''}
\providecommand{\smfandname}{et}
\providecommand{\smfedsname}{\'eds.}
\providecommand{\smfedname}{\'ed.}
\providecommand{\smfmastersthesisname}{M\'emoire}
\providecommand{\smfphdthesisname}{Th\`ese}
\begin{thebibliography}{10}

\bibitem{Bouchitte}
{\scshape G.~Bouchitt{\'e}, C.~Jimenez {\normalfont \smfandname} M.~Rajesh} --
  {\og Asymptotique d'un probl\`eme de positionnement optimal\fg}, \emph{C. R.
  Math. Acad. Sci. Paris} \textbf{335} (2002), no.~10, p.~853--858.

\bibitem{Brancolini}
{\scshape A.~Brancolini, G.~Buttazzo, F.~Santambrogio {\normalfont \smfandname}
  E.~Stepanov} -- {\og Long-term planning versus short-term planning in the
  asymptotical location problem\fg}, \emph{ESAIM Control Optim. Calc. Var.}
  \textbf{15} (2009), no.~3, p.~509--524.

\bibitem{Champion-De_Pascale-Petri}
{\scshape T.~Champion, L.~De~Pascale {\normalfont \smfandname} P.~Juutinen} --
  {\og The {$\infty$}-{W}asserstein distance: local solutions and existence of
  optimal transport maps\fg}, \emph{SIAM J. Math. Anal.} \textbf{40} (2008),
  no.~1, p.~1--20.

\bibitem{Dobric-Yukich}
{\scshape V.~Dobri{\'c} {\normalfont \smfandname} J.~E. Yukich} -- {\og
  Asymptotics for transportation cost in high dimensions\fg}, \emph{J. Theoret.
  Probab.} \textbf{8} (1995), no.~1, p.~97--118.

\bibitem{Du-Faber-Gunzburger}
{\scshape Q.~Du, V.~Faber {\normalfont \smfandname} M.~Gunzburger} -- {\og
  Centroidal {V}oronoi tessellations: applications and algorithms\fg},
  \emph{SIAM Rev.} \textbf{41} (1999), no.~4, p.~637--676 (electronic).

\bibitem{Du-Wang}
{\scshape Q.~Du {\normalfont \smfandname} D.~Wang} -- {\og The optimal
  centroidal {V}oronoi tessellations and the {G}ersho's conjecture in the
  three-dimensional space\fg}, \emph{Comput. Math. Appl.} \textbf{49} (2005),
  no.~9-10, p.~1355--1373.

\bibitem{Falconer}
{\scshape K.~J. Falconer} -- \emph{The geometry of fractal sets}, Cambridge
  Tracts in Mathematics, vol.~85, Cambridge University Press, Cambridge, 1986.

\bibitem{FejesToth}
{\scshape L.~Fejes~T{\'o}th} -- {\og Sur la repr\'esentation d'une population
  infinie par un nombre fini d'\'el\'ements\fg}, \emph{Acta. Math. Acad. Sci.
  Hungar} \textbf{10} (1959), p.~299--304 (unbound insert).

\bibitem{FejesToth2}
{\scshape L.~Fejes~T{\'o}th} -- \emph{Lagerungen in der {E}bene auf der {K}ugel
  und im {R}aum}, Springer-Verlag, Berlin, 1972, Zweite verbesserte und
  erweiterte Auflage, Die Grundlehren der mathematischen Wissenschaften, Band
  65.

\bibitem{Graf-Luschgy}
{\scshape S.~Graf {\normalfont \smfandname} H.~Luschgy} -- \emph{Foundations of
  quantization for probability distributions}, Lecture Notes in Mathematics,
  vol. 1730, Springer-Verlag, Berlin, 2000.

\bibitem{Heinonen}
{\scshape J.~Heinonen} -- \emph{Lectures on analysis on metric spaces},
  Universitext, Springer-Verlag, New York, 2001.

\bibitem{Horowitz-Karandikar}
{\scshape J.~Horowitz {\normalfont \smfandname} R.~L. Karandikar} -- {\og Mean
  rates of convergence of empirical measures in the {W}asserstein metric\fg},
  \emph{J. Comput. Appl. Math.} \textbf{55} (1994), no.~3, p.~261--273.

\bibitem{Hutchinson}
{\scshape J.~E. Hutchinson} -- {\og Fractals and self-similarity\fg},
  \emph{Indiana Univ. Math. J.} \textbf{30} (1981), no.~5, p.~713--747.

\bibitem{Morgan-Bolton}
{\scshape F.~Morgan {\normalfont \smfandname} R.~Bolton} -- {\og Hexagonal
  economic regions solve the location problem\fg}, \emph{Amer. Math. Monthly}
  \textbf{109} (2002), no.~2, p.~165--172.

\bibitem{Mosconi-Tilli}
{\scshape S.~J.~N. Mosconi {\normalfont \smfandname} P.~Tilli} -- {\og
  {$\Gamma$}-convergence for the irrigation problem\fg}, \emph{J. Convex Anal.}
  \textbf{12} (2005), no.~1, p.~145--158.

\bibitem{Newman}
{\scshape D.~J. Newman} -- {\og The hexagon theorem\fg}, \emph{IEEE Trans.
  Inform. Theory} \textbf{28} (1982), no.~2, p.~137--139.

\bibitem{Villani}
{\scshape C.~Villani} -- \emph{Topics in optimal transportation}, Graduate
  Studies in Mathematics, vol.~58, American Mathematical Society, Providence,
  RI, 2003.

\end{thebibliography}

\end{document}